\documentclass{elsart3-1}

\usepackage{latexsym}
\usepackage{amsmath}
\usepackage{amsfonts}
\usepackage{amssymb}
\usepackage{graphpap}
\usepackage{epsfig}
\usepackage{amscd}
\usepackage{enumerate}
\usepackage{eucal} 

\newcommand{\RR}{{\mathbb{R}}}
\newcommand{\NN}{{\mathbb{N}}}
\newcommand{\ZZ}{{\mathbb{Z}}}

\newcommand{\eps}{\varepsilon}

\newcommand{\ep}{\hfill$\Box$ \vskip 0.08in}

\usepackage[english,francais]{babel}

\newtheorem{theorem}{Theorem}[section]

\newtheorem{e-proposition}[theorem]{Proposition}

\newtheorem{e-definition}[theorem]{Definition\rm}

\newtheorem{theoreme}{Th\'eor\`eme}[section]

\newtheorem{proposition}[theoreme]{Proposition}

\renewcommand{\theequation}{\arabic{equation}}
\def\og{\leavevmode\raise.3ex\hbox{$\scriptscriptstyle\langle\!\langle$~}}
\def\fg{\leavevmode\raise.3ex\hbox{~$\!\scriptscriptstyle\,\rangle\!\rangle$}}

\journal{}
\begin{document}
\centerline{}
\begin{frontmatter}


\selectlanguage{english}
\title{Square function and Riesz transform  in non-integer dimensions}


\selectlanguage{english}
\author[SM]{Svitlana Mayboroda},
\ead{svitlana@math.purdue.edu}
\author[AV]{Alexander Volberg}
\ead{volberg@math.msu.edu, A.Volberg@ed.ac.uk}

\address[SM]{Department of Mathematics, Purdue University, 150 N. University Street, West Lafayette, IN 47907-2067, USA}
\address[AV]{Department of Mathematics, Michigan State University, East Lansing, MI 48824, USA}


\medskip
\begin{center}
\end{center}

\begin{abstract}
\selectlanguage{english}
Following a recent paper \cite{RdVTolsa} we show that the finiteness of square function associated with the Riesz transforms with respect to Hausdorff measure $H^s$ implies that $s$ is integer. 

\vskip 0.5\baselineskip

\selectlanguage{francais}

\end{abstract}
\end{frontmatter}


\selectlanguage{english}
\section{Introduction}
\label{}


\font\tretamajka=cmbsy10 at 12pt

For a Borel mesure $\mu$ in $\RR^m$ and $s\in (0,m]$ the $s$-Riesz transform of $\mu$ is defined as
\begin{equation}\label{eq1.1}
R^s \mu(x):=\int\frac{x-y}{|x-y|^{s+1}}\,d\mu(y), \qquad x\not\in {\rm supp}\,\mu,
\end{equation}

\noindent and the truncated Riesz transform is given by
\begin{equation}\label{eq1.2}
R^s_{\eps}\, \mu(x):=\int_{|x-y|>\eps}\frac{x-y}{|x-y|^{s+1}}\,d\mu(y), \quad R^s_{\eps,\eta}\, \mu(x):=R^s_{\eta}\, \mu(x)-R^s_{\eps}\, \mu(x),
\end{equation}

\noindent where $x\in\RR^m,$ $\eta>\eps>0$.

Further, recall that the upper and lower $s$-dimensional densities of $\mu$ at $x$ are given by  
\begin{equation}\label{eq1.2.1}
\theta^{s,\ast}_{\mu}(x):=\limsup_{r\to 0}\frac{\mu(B(x,r))}{r^s} \quad\mbox{and}\quad \theta^{s}_{\mu,\ast}(x):=\liminf_{r\to 0}\frac{\mu(B(x,r))}{r^s},
\end{equation}

\noindent respectively, where $B(x,r)$ is the ball of radius $r>0$ centered at $x\in \RR^m$. 

It has been proved in \cite{Pr1} and \cite{Pr2} that whenever $0\leq s\leq 1$ and $\mu$ is a finite Radon measure with $0<\theta^{s,\ast}_{\mu}(x)<\infty$, for $\mu - {\rm a.e.}\,x\in\RR^m$,
the condition 
\begin{equation}\label{eq1.2.2}
\sup_{\eps>0} |R_{\eps}^s\,\mu(x)|<\infty\qquad \mu - {\rm a.e.}\,x\in\RR^m,
\end{equation}

\noindent implies that $s\in\ZZ$. Moreover, an analogous result has been obtained in \cite{Vi} for all $0\leq s\leq m$ under a stronger assumption that  $0<\theta^{s}_{\mu,\ast}(x)\leq \theta^{s,\ast}_{\mu}(x) <\infty$. However, neither the curvature methods of \cite{Pr1}, \cite{Pr2}, nor the tangent measure techniques in \cite{Vi} could be applied to establish that \eqref{eq1.2.2} implies $s\in\ZZ$ for all $0\leq s\leq m$ assuming only $0<\theta^{s,\ast}_{\mu}(x)<\infty$.

In  \cite{RdVTolsa} the authors proved that the latter statement holds if the condition \eqref{eq1.2.2} is substituted by the existence of  the principal value $\lim_{\eps\to 0} R_{\eps}^s\,\mu(x)$, $\mu - {\rm a.e.}\,x\in\RR^m$. In the present work we refine the techniques of  \cite{RdVTolsa} and establish the following result. 

\begin{theorem}\label{tSqF} Let $\mu$ be a finite Radon measure in $\RR^m$ with 
\begin{equation}\label{eq1.2.3}
0<\theta^{s,\ast}_{\mu}(x)<\infty\quad\mbox{for}\quad \mu - {\rm a.e.}\,x\in\RR^m.
\end{equation}

\noindent Furthermore, assume that  for some $s\in (0,m]$ the square function
\begin{equation}\label{eq1.3}
S^s\mu(x):=\Bigg(\int_0^\infty \left|R^s_{t,2t}\, \mu(x)\right|^2\,\frac{dt}{t}\Bigg)^{1/2}<\infty, \qquad \mu - {\rm a.e.}\,x\in\RR^m.
\end{equation}

\noindent Then $s\in\ZZ$.
\end{theorem}


In fact, we also prove the following closely related result which is a strengthening of  the main results in \cite{RdVTolsa}.\begin{theorem}\label{tLim} Let $\mu$ be a finite Radon measure in $\RR^m$  satisfying \eqref{eq1.2.3}. Furthermore, assume that  for some $s\in (0,m]$ we have
\begin{equation}\label{eq1.3.1}
\lim_{\eps\to 0}R^s_{\eps,2\eps}\,\mu(x)=0 \qquad \mu - {\rm a.e.}\,x\in\RR^m.
\end{equation}

\noindent Then $s\in\ZZ$.
\end{theorem}

This circle of problems goes back, in particular, to the work of David and Semmes \cite{DS1}, \cite{DS2}, where the authors showed, under certain assumptions on the measure $\mu$, that the $L^2$ boundedness of a large class of  singular integral operators implies that $s$ is an integer and $\mu$ is uniformly rectifiable, that is, the support of $\mu$ contains ``large pieces of Lipschitz graphs" -- see \cite{DS1}, \cite{DS2} for details. The ultimate goal, which seems to be out of reach at the moment, is to prove that a similar conclusion holds purely under the assumption that the Riesz transform is bounded in $L^2$, i.e., that the Riesz transform alone encodes the geometric information about the underlying measure. The achievements in \cite{NTV-Acta}, \cite{VBook}, \cite{NTV-IMRNWeakType} showed that the $L^2$-boundedness of the Riesz transform, suitably interpreted, is almost equivalent to the condition \eqref{eq1.2.2}. However, under the assumption \eqref{eq1.2.2} the problem seems to be just as challenging.
In both cases the question has only been resolved for $s=1$ (\cite{MMV}, \cite{To2}, \cite{To10}), by the methods involving curvature of measures. 

In this vein, we would like to point out that by Khinchin's inequality \eqref{eq1.3} can be viewed {\it almost} as a condition
\begin{equation}\label{eq1.4}
{\mathbb{E}}\,\Bigg|\sum_{k\in\ZZ} \eps_k R^s_{2^{-k},2^{-k+1}}\, \mu(x)\Bigg|<\infty \qquad \mu - {\rm a.e.}\,x\in\RR^m,
\end{equation}

\noindent where $\eps_k$ are independent random variables taking the values $-1$ and $1$ with probability $1/2$ each. Therefore,  in order to guarantee $s\in\ZZ$, it is sufficient to assume only that the singular integrals of the type $\sum_{k=0}^{\infty} \eps_k R^s_{2^{-k},2^{-k+1}}\, \mu(x)$ are uniformly bounded.

Finally, the $s$ dimensional Hausdorff measure $H^s$ of a set $E$ with $0<H^s(E)<\infty$ satisfies the condition \eqref{eq1.2.3}, and hence, the results of Theorems~\ref{tSqF} and \ref{tLim} remain valid in this context.

\section{Preliminary estimates}

Our proof largely relies on the estimates for a slightly modified version of the Riesz transform that were obtained in \cite{RdVTolsa}. To be precise, let us consider the operator
\begin{equation}\label{eq2.1}
R^{s,\varphi}_{\eps}\, \mu(x):=\int\varphi\left(\frac{|x-y|^2}{\eps^2}\right)\,\frac{x-y}{|x-y|^{s+1}}\,d\mu(y), \qquad \eps>0,
\end{equation}

\noindent where $\varphi=\varphi_\rho$ is a $C^2$ function depending on the parameter $\rho\in(0,1/2)$, to be determined below, with ${\rm supp}\,\varphi\subset [0,1+\rho+2\rho^2]$ and such that
\begin{enumerate}
\item $\varphi(r)=r^{\frac{s+1}{2}}$ for $0\leq r \leq 1$,
\qquad $\varphi(r)=-\frac r\rho+1+\rho+\frac 1\rho$ for $1+\rho^2\leq r\leq 1+\rho^2+\rho$,
\item $|\varphi(r)|\leq C,$ $|\varphi'(r)|\leq 1/\rho$, $|\varphi''(r)|\leq C(\rho)$ for all $r>0$.
\end{enumerate} 

Analogously to \cite{RdVTolsa}, we start with a set 
\begin{eqnarray}\label{eq2.2}
\nonumber F_\delta &:=& \{x\in \RR^m:\,\mu (B(x,r))/r^s\leq 2\theta_\mu^{s,\ast}(x)\,\mbox{ for }\,r\leq r_0,\quad \theta_\mu^{s,\ast}(x)\leq C_0,\\[4pt]
&&\quad \mbox{and } |R^{s,\varphi}_{\eps}\,\mu(x)-R^{s, \varphi}_{2\eps}\,\mu(x)|\leq \delta \,\mbox{ for all }\,0<\eps<\eps_0\},
\end{eqnarray}

\noindent where $0<\delta<1$ and $C_0, r_0,\eps_0$ are some positive constants.  

One can see that both the condition \eqref{eq1.3} and \eqref{eq1.3.1} imply that 
\begin{equation}\label{eq2.4}
\lim_{\eps\to 0}|R^{s,\varphi}_{\eps}\,\mu(x)-R^{s,\varphi}_{2\eps}\,\mu(x)|=0 \qquad \mu - {\rm a.e.}\,x\in\RR^n.
\end{equation}
\noindent Indeed, 
\begin{equation}\label{eq2.3}
R^{s,\varphi}_{\eps}\, \mu(x)=\int\int_{0<t<\frac{|x-y|^2}{\eps^2}}\varphi'(t)\,dt\,\frac{x-y}{|x-y|^{s+1}}\,d\mu(y)=\int_0^{1+\rho+2\rho^2} \varphi'(t)R_{\eps\sqrt t}^{s}\,\mu(x)\,dt, 
\end{equation}

\noindent so that  $\eqref{eq1.3.1}$ directly gives \eqref{eq2.4}. Furthermore,  $\eqref{eq2.3}$ entails that
\begin{equation}\label{eq2.3.1}\nonumber
\left|R^{s,\varphi}_{\eps}\, \mu(x)-R^{s,\varphi}_{2\eps}\, \mu(x)\right|\leq C(\rho) \int_0^{1+\rho+2\rho^2} |R_{\eps\sqrt t,2\eps\sqrt t}^{s}\,\mu(x)|\,dt\leq C(\rho)\left(\int_0^{\eps\sqrt{1+\rho+2\rho^2}}|R_{u,2u}^{s}\,\mu(x)|^2\,\frac{du}{u}\right)^{1/2},
\end{equation}

\noindent where we used the change of variables $u:=\eps\sqrt t$ and H\"older's inequality. Hence, \eqref{eq1.3} also leads to \eqref{eq2.4}.

 Therefore, for sufficiently small $\eps_0$ and $r_0$ and sufficiently large $C_0$ the set $F_\delta$ has $\mu(F_{\delta})>0$. 
Note that $\mu (B(x,r))\leq Mr^s$ for all $x\in F_\delta$, $r>0$ and $M=\max\{2C_0,\mu(\RR^m)/r_0^s\}$. 

Let $\theta^s(x,r)$ denote  the average $s$-dimensional density of the ball $B(x,r)$, $x\in\RR^m$, $r>0$, that is, 
$\theta^s(x,r):=\mu(B(x,r))/r^s$. We start with the following estimates.

\begin{proposition}\cite{RdVTolsa} \label{p2.1}
Assume that for some $C'>0$, $r>0$ and $x_0\in\RR^m$ we have 
$\mu(B(x_0,r))\geq C'r^s$, and denote by $n$ the biggest integer strictly smaller than $s$. Then for a sufficiently small $\rho$ (depending on $s$ only) and any $\tau\in (0,1/20)$ there exists a constant $\omega_0=C(C',M,\rho)\tau^{-s\,\frac {1}{\log_4(1+\rho^2/4)}}$, there exists $\eps\in\left(\frac r\tau, \omega_0\,\frac r\tau\right]$ and a set of points $y_0,...,y_{n+1}\in B(x_0,r)\cap F_{\delta}$ such that 
\begin{equation}\label{eq2.9}
\theta^s(y_0,4\eps)\leq C(\rho)\,\theta^s(y_0,\eps),\qquad \theta^s(y_0,\eps)\geq C'\tau^s/2,
\end{equation}

\noindent and 
\begin{equation}\label{eq2.7}
\sum_{j=1}^{n+1}|R^{s,\varphi}_{\eps}\,\mu(y_j)-R^{s,\varphi}_{\eps}\,\mu(y_0)|+\theta^s(y_0, 3\eps)\frac{r^2}{\eps^2}\geq C(C',M,s)(n+1-s)\,r\,\frac{\theta^s(y_0,\eps)}{\eps}.
\end{equation}
\end{proposition}

\section{The proof of the main result}

 
We will argue by contradiction. We initially assume that $s\not\in \ZZ$, and then show that the estimate in \eqref{eq2.7} is accompanied by the corresponding bound from above in terms of $\delta$, $r$, $\tau$ and $\eps$. Ultimately, choosing $r$, $\delta$, $\tau$ sufficiently small leads to a contradiction. Observe that in the case $s\in\ZZ$ the lower bound in \eqref{eq2.7} is degenerate, and hence, such an argument could not be constructed.

Set 
\begin{equation}\label{eq3.1}
\delta:=\frac{\tau^{s+2}}{\omega_0}=C(M,\rho)\, \tau^{s+2+s\,\frac {1}{\log_4(1+\rho^2/4)}},
\end{equation}

\noindent where the constant $C(M,\rho)$ is equal to the reciprocal of $C(C',M,\rho)$ from the definition of $\omega_0$ corresponding to $C'=1/2$. Note that $M$ depends on $r_0$ and $C_0$ in the definition of $F_\delta$, however, the choice of $r_0$ and $C_0$ is determined solely by the properties of $\mu$ and can be made independent of $\delta$.

Going further, fix $\eps_0$ and take 
\begin{equation}\label{eq3.2}
r<\eps_0\delta=\eps_0\,\frac{\tau^{s+2}}{\omega_0}=C(M,\rho)\,\eps_0\, \tau^{s+2+s\,\frac {1}{\log_4(1+\rho^2/4)}}
\quad \mbox{such that} \quad \mu(B(x_0,r)\cap F_\delta)\geq r^s/2.\end{equation}

Now that $r$ and $\delta$ are fixed, we invoke Proposition~\ref{p2.1}, find the points $y_0,..., y_{n+1}$  and choose $\eps\in \left(\frac r\tau,\omega_0\,\frac r\tau\right]$ such that \eqref{eq2.9} and \eqref{eq2.7} are satisfied. However,  for every $x,z\in B(x_0, r)\cap F_{\delta}$, $x\in \RR^m$ we have
\begin{equation}\label{eq3.4}
|R^{s,\varphi}_{\eps}\,\mu(x)-R^{s,\varphi}_{\eps}\,\mu(z)|\leq C(\rho)M \delta + C \delta \log \frac {r}{\delta\eps},
\end{equation}

\noindent whenever  $2\eps<\frac r\delta<\eps_0$. Indeed, a direct calculation shows that for $\eta\in\left[\frac{r}{2\delta},\frac{r}{\delta}\right]$
\begin{equation}\label{eq3.6}
|R^{s,\varphi}_{\eta}\,\mu(x)-R^{s,\varphi}_{\eta}\,\mu(z)|\leq C(\rho)\left(r/\delta\right)^{-s-1}\, |z-x| \,\mu(B(x_0, 4r/\delta))\leq C(\rho)M \delta.
\end{equation}

\noindent Then we can choose $\eta\in\left[\frac{r}{2\delta},\frac{r}{\delta}\right]$ such that $\eta=2^k\eps$ for some $k\in\NN$, so that
\begin{eqnarray}\label{eq3.7}
&&|R^{s,\varphi}_{\eps}\,\mu(x)-R^{s,\varphi}_{\eps}\,\mu(z)|\nonumber\\[4pt]\nonumber &&\qquad\leq |R^{s,\varphi}_{\eps}\,\mu(x)-R^{s,\varphi}_{\eta}\,\mu(x)|+
|R^{s,\varphi}_{\eta}\,\mu(x)-R^{s,\varphi}_{\eta}\,\mu(z)|+|R^{s,\varphi}_{\eta}\,\mu(z)-R^{s,\varphi}_{\eps}\,\mu(z)| \\[4pt] &&\qquad
\leq C(\rho)M \delta + C \sup_{x\in F_{\delta}}\sup_{1\leq i\leq k}|R^{s,\varphi}_{2^{i}\eps}\,\mu(x)-R^{s,\varphi}_{2^{i-1}\eps}\,\mu(x)| \log \frac {r}{\delta\eps}\leq C(\rho)M \delta + C \delta \log \frac {r}{\delta\eps}.
\end{eqnarray}

Therefore, \eqref{eq2.7} is complemented by the estimate
\begin{equation}\label{eq3.5}
\sum_{j=1}^{n+1}|R^{s,\varphi}_{\eps}\,\mu(y_j)-R^{s,\varphi}_{\eps}\,\mu(y_0)|+\theta^s(y_0, 3\eps)\frac{r^2}{\eps^2}\leq 
C(M,\rho)\left(\delta+\delta \log \frac{r}{\eps\delta}+ \theta^s(y_0, \eps)\frac{r^2}{\eps^2}\right),\end{equation}

\noindent where we used \eqref{eq3.4} and \eqref{eq2.9}. Now combining \eqref{eq2.7} with \eqref{eq3.5} and dividing both sides by $r/\eps$ we arrive at the estimate
\begin{equation}\label{eq3.8}
\theta^s(y_0,\eps)\leq C(M,s,\rho)\left(\frac{\delta\eps}{r}+\frac{\delta\eps}{r}\log \frac{r}{\eps\delta}+\theta^s(y_0, \eps)\frac{r}{\eps}\right).
\end{equation}

\noindent According to our choice of $\delta$ and $r$, 
\begin{equation}\label{eq3.9}
\frac{\delta\eps}{r}\leq \frac{\tau^{s+2}}{\omega_0}\,\frac{\omega_0}{\tau}=\tau^{s+1}\leq C\tau \,\theta^s(y_0, \eps) \quad \mbox{and}\quad \frac{r}{\eps} \leq \tau.
\end{equation}

 Now \eqref{eq3.8} and \eqref{eq3.9} give the bound
\begin{equation}\label{eq3.10}
\theta^s(y_0,\eps)\leq C(M,s,\rho)\left(\tau +\tau^{1-\alpha}\right)\theta^s(y_0, \eps), \qquad \forall\,\alpha>0,
\end{equation}

\noindent which for $\tau>0$ sufficiently small leads to a contradiction. \ep

\bigskip

\end{document}